
\magnification=\magstep1
\documentstyle{amsppt}
\pagewidth{6.5truein}
\pageheight{9.0truein}
\NoBlackBoxes
\TagsOnRight
\centerline {\bf Lines on Calabi Yau complete intersections, mirror
 symmetry,}
\centerline {\bf and Picard Fuchs equations}
\bigskip
\centerline{\bf A. Libgober\footnote{Supported by NSF grant DMS--910798.}
 and J. Teitelbaum\footnote{Supported by NSF grant DMS--9204265 and by a Sloan
 Research Fellowship.}}
\centerline {Department of Mathematics}
\centerline {University of Illinois at Chicago}
\centerline {P.O.B. 4348, Chicago, Ill, 60680}
\centerline{Internet: u11377\@uicvm.uic.edu\quad jeremy\@math.uic.edu}
\bigskip
\par {\bf Introduction and statement of the result.}
It was suggested (cf. [COGP], [GP])  that (in some circumstances)
if $V$ is a Calabi--Yau threefold
 then one can  relate to $V$ a family $W(V)_t$, $t \in {\bold C}$
  of Calabi--Yau manifolds which are "mirrors" of $V$, such that
 one has the following relation between their Euler characteristics:
 $\chi (V)=-\chi (W(V)_t)$.  One
  of the properties of this correspondence should be the
 following: the coefficients of the expansion of certain integrals
 attached to $W(V)_t$ (so called  Yukawa couplings) relative to an
 appropriately chosen parameter are  integers from which one may calculate the
 numbers
 $r_d$ of rational curves  of  degree $d$ on
 a generic Calabi--Yau manifold which is a deformation of $V$.
   This was verified in [COGP] and [M1],[M2]
  in the case when $V$ is the
 quintic hypersurface in ${\bold{CP}}^4$ for rational curves of a
 small degree.
 Other authors ([CL],[S]) have  suggested a large list of mirrors of
 hypersurfaces in
 weighted projective spaces. The purpose of this note is
to verify the above predictions for the  remaining types of  Calabi--Yau
{\it  complete intersections}
 in complex projective space when $d=1$ i.e. the case of lines.
\par  A Calabi--Yau threefold $W$ is a Kahler manifold such that $\dim W=3$,
the canonical bundle of $W$ is trivial, and the Hodge numbers satisfy
$h^{1,0}=h^{2,0}=0$. Let $W_{t}$ be a family of such manifolds and let
 $\omega_t$ be a
 family of   holomorphic 3-forms on $W_t$
   (unique up to constant for each $t$ because $h^{3,0}(W_t)=1$).
 According to Griffiths transversality ([G]):
 $$\kappa_t(k)=\int_W \omega \wedge
 {{d^k \omega_t} \over {d t^k}} \eqno (1)$$
 is equal to zero for $k \le 2$. Let
$\kappa_{ttt}= \kappa_t(3)$.
We assume that
 the monodromy $T$ about $t=\infty$ acting on
  $H_3 (W_t,{\bold Z})$  is maximally unipotent i.e. that $(T-I)^3 \ne 0$
 and $(T-I)^4=0$. If this is the case then ([M1],[M2]) for $N=logT$
one has $dim (Im N^3)  \otimes {\bold C}=1$ ( as a consequence of
 $h^{3,0}=1$). Let $\gamma_1,\gamma_0 \in H_3 (W_t,{\bold Z})$
  be a basis of $(Im  N^2)
\otimes {\bold C}$ such that $\gamma_0 \in (Im N^3)$ is an indivisible
element and $\gamma_1=1/\lambda N^2 \tilde \gamma_1$ where $\gamma_1$
 is indivisible and the intersection index of $\tilde \gamma_1$ and
 $\gamma_0$ is 1. Let $m$ be defined from the relation
 $N \gamma_1= m \cdot \gamma_0$ and let $$s={{{1 \over m}
 \int_ {\gamma_1} \omega } \over {\int _{\gamma_0} \omega }},
 \   q=e^{2 \pi i s}. \eqno (2)$$
 Then $q$ is indepenedent of a choice of the basis
 $\gamma_0, \gamma_1$ and the form $\omega$
  up to root of unity of degree $\vert m \vert$ (cf. [M1]).
\par In this paper, we use an {\it a priori}  different normalization for the
 parameter $s$
determined by specifying the asymptotic behavior of $s$ as $t\to\infty$.  This
 normalization
is described in (18); it is analogous to that exploited in [COGP] and [M2].
\par  Let $V_{\lambda}$ be given in ${\bold{CP}}^5$ by
 $$ Q_1=x_1^3+x_2^3+x_3^3-3 \lambda x_4 x_5 x_6=0$$
 $$Q_2=x_4^3+x_5^3+x_6^3-3 \lambda x_1 x_2 x_3 =0 \eqno (3)$$
 This is a complete intersection which is a Calabi--Yau threefold
 for generic $\lambda$ . Let $G_{81} \subset PGL(5,{\bold C})$
  be the subgroup (of order $81$)
 of transformations $g_{\alpha,\beta,\delta,\epsilon,\mu}$
 where $\alpha,\beta,\delta,\epsilon \in {\bold Z} \  \bmod 3$ $\mu \in
 {\bold Z}_9$ and $3 \cdot \mu =\alpha +\beta=\delta+ \epsilon \ \bmod 3$.
 These transformations act as:
 $$
g_{\alpha,\beta,\delta,\epsilon,\mu}:
 (x_1,x_2,x_3,x_4,x_5,x_6) \rightarrow $$ $$\rightarrow
 (\zeta_3^\alpha \cdot \zeta _9
 ^\mu \cdot x_1, \zeta_3^{\beta} \cdot \zeta_9^{\mu} \cdot x_2,
 \zeta_9^{\mu} \cdot x_3, \zeta _3^{-\delta} \cdot  \zeta_9 ^{-\mu}
\cdot x_4, \zeta _3 ^{-\epsilon} \cdot \zeta_9 ^{-\mu} \cdot x_5,
 \zeta_9 ^{-\mu} \cdot x_6) \eqno (4)$$
 and preserve both hypersurfaces $\tilde Q_i$ given by the equations
  $Q_i=0$ $(i=1,2)$.
\bigskip {\bf Theorem}. {\it The resolution of
 singularities $W(V_\lambda)$
 of the quotient of $V_{\lambda}$  by the action of $G_{81}$
 which is a Calabi--Yau manifold  satisfies: $\chi (V_\lambda)=
 -\chi (W(V_\lambda))$. The monodromy of $W(V_\lambda)$ about
 infinity is maximally unipotent.
 For $q$ defined for the family
 $W(V_\lambda)$ by the asymptotic normalization (18),
the coefficient of  $q$ in the $q$-expansion of
 $\kappa_{sss}$ is equal to the number of
 lines on a generic non singular complete intersection of two
 cubic hypersurfaces in ${\bold{CP}}^5$.}
\bigskip {\bf A calculation of the Euler characteristic.}
 The statement on Euler characteristic  (as well as the statement on
   the number of lines)  are verified
 by direct calculation of the quantities involved. The total Chern class
  $c=1+c_1 +c_2  +c_3 \in H^* (V_{\lambda}, {\bold Z})$, of the tangent
 bundle of $V_{\lambda}$  satisfies
 $c \cdot (1+3 \cdot h)^2= (1+h)^6$ where $h$ is the generator
 $H^2 (V_{\lambda},{\bold Z})$ (here $(1+3 \cdot h)^2$ and $(1+h)^6$
 respectively the total Chern class of the normal bundle to
 $V_\lambda$  in ${\bold{CP}}^5$ and the pulback on $V_{\lambda}$
  of the total Chern class of ${\bold{CP}}^5$) . The Euler charcateristic
 of $V_{\lambda}$ is  $c_3$ evaluated on its fundamental class
which (using the fact that $h^3$ evaluated on the fundamental class is
 $9$) gives  $\chi (V_\lambda)=-144$.
\par On the other hand according to the ``physicist's formula'' (cf. [DHVW])
 or rather to its reformulation due to Hirzebruch and Hofer (cf. [HH])
 the Euler characteristic of a  Calabi--Yau resolution of the quotient
 $V_{\lambda}/ G_{81}$ can be found as
 $$\Sigma_{[g]} \chi (V_{\lambda}^g/C(g)) \eqno (5)$$
 where  the summation is over all conjugacy classes $[g]$ of elements of $G$,
 $C(g)$ denotes the centralizer of $g$ and $X^g$ is the fixed point set
 of an element $g$. Because $G_{81}$ is abelian the formula reduces to
 $\Sigma_g \chi (V_{\lambda}^g/G_{81})$
 where the summation is over all elements of
  the group.
  There are $6$ curves $C_{i,j}$ having non--trivial stabilizer
 corresponding to the vanishing of two variables in either of the two sets
 $(x_1,x_2,x_3)$ or $(x_4,x_5,x_6)$. The Euler
 characteristic of such a curve, which is a complete intersection
 of two cubic surfaces in ${\bold P}^3$, is $-18$.
 The stabilizer of each
 curve contains
  $3$ elements since for each curve there are $2$ elements which
have this curve as the fixed point set. Hence the number of elements
which have one dimensional fixed point set is 12. The Euler
 characteristic of the quotient of the one dimensional fixed point
 set  is $2$. The zero dimensional fixed point sets $D_{i,j,k}$
 on a curve $C_{i,j}$
 ($i,j$ are in the same group of variables, and $k$ in another)
  are obtained
 by equating to zero a variable in another group. The stabilizer
 of such zero dimensional fixed point set  has order $27$. Each
 zero dimensional fixed point set $D_{i,j,k}$  belongs to $3$
 curves $C_{i,j}$. Hence the number of elements stabilizing $D_{i,j,k}$
  is $27-3 \times 2-1=20$. The number of zero dimensional fixed point
  sets $D_{i,j,k}$ is $6$  and each element with zero dimensional fixed
point set stabilizes $2$ sets $D_{i,j,k}$. Hence  the number of elements
 with zero dimensional stabilizer is $20 \times 6/2=60$ and each
 such element stablizes $6$ points.  The quotient
 of a zero dimensional fixed point set of an element by the group has the
 Euler characteristic equal to $2$. The contribution in (5) from
the identity element is $$\chi (V_{\lambda}/G_{81})={1 \over
 {{\vert G_{81}
 \vert}}} \Sigma _g \chi (V_{\lambda}^g)$$ which is equal to
  $1/81 (-144+60 \times 6 + 12 \times (-18))=0$ .
 Hence using (5) the Euler
characteristic of a Calabi--Yau resolution is $0+2 \times 60+2 \times
 12=+144$.
\bigskip {\bf A method for constructing the Picard Fuchs equations.}
 To find the Picard Fuchs equations for the periods of $W(V_\lambda)$
we shall extend
  to complete intersections  the Griffiths
   description of cohomology classes of
  hypersurfaces using meromorphic forms on the ambient space.
  Let $T(\tilde  Q_1 \cap \tilde  Q_2)$
    be a small tubular  neighbourhood of $\tilde Q_1 \cap \tilde Q_2$ in
 ${\bold{CP}}^5$ and $\partial (T(\tilde Q_1 \cap \tilde Q_2))$ be the
 boundary of $T(\tilde  Q_1 \cap \tilde Q_2)$. Then

$$
H^3 (\tilde Q_1 \cap \tilde Q_2)
		  = H^3 (\tilde Q_1 \cap \tilde Q_2)^*
		 = H^7 (T(\tilde Q_1 \cap \tilde Q_2,\partial T(\tilde Q_1 \tilde Q_2))
		 = H^7 ({\bold{CP}}^5,{\bold{CP}}^5-T(\tilde Q_1 \cap \tilde Q_2))
$$
 (use Poincare duality, retraction combined with  Lefschetz duality,
  and  excision). The latter group is isomorphic to $H^6 ({\bold{CP}}^5-
 \tilde Q_1 \cap \tilde Q_2)$ as follows from the exact sequence of the pair.
 The Mayer Vietoris  sequence combined with these isomorphisms gives
 the  identification:
 $$H^5 ({\bold{CP}}^5- (\tilde Q_1 \cup \tilde Q_2)/ Im \ (
 H^5 ({\bold{CP}}^5-\tilde Q_1) \oplus H^5({\bold{CP}}^5- \tilde Q_2))=
    H^3 (\tilde Q_1 \cap \tilde Q_2) \eqno (7)$$
An alternative description of this isomorphism can be obtained by
interpreting
  a meromorphic 5-form on ${\bold{CP}}^5$ having poles along $\tilde Q_1
 \cup \tilde Q_2$  as a functional on $H_3 (\tilde
 Q_1 \cap \tilde Q_2)$ which is given by assigning to a 3-cycle
  $\gamma$ representing a homology class in the latter group the integral
 over a 5-cycle in ${\bold{CP}}^5-(\tilde Q_1 \cup \tilde Q_2)$; This 5-cycle
 is the restriction to $\gamma$ of a torus fibration on which
 $T(\tilde Q_1 \cap \tilde Q_2) - (\tilde Q_1 \cup \tilde Q_2)
 \cap T(\tilde Q_1 \cap \tilde Q_2)$ retracts as a consequence of
 the non--singularity of $\tilde Q_1 \cap \tilde Q_2$. Moreover in the
 isomorphism (7) the filtration by the total order of the pole
 corresponds  to the Hodge filtration on $H^3 (\tilde Q_1 \cap
 \tilde Q_2)$ (details of this will appear elsewhere).
The residues of the meromorphic 5-forms which are $G_{81}$-invariant
  give the forms  on $V_\lambda$ which
descend to $V_\lambda/G_{81}$; The  pull--back of these forms, which
 give a basis of $H^3 (W(V_\lambda)$, are
   $${(x_1x_2x_3)^{i-1}(x_4x_5x_6)^{n-i-1}}\Omega \over {Q_1^iQ_2^{n-i}}
    \eqno (8)$$
  where $n=2,3,4,5$ and $\Omega$ is the Euler form:
$$
\Omega=\sum (-1)^{i}x_{i} dx_{1}\wedge\cdots\wedge \hat{dx_{i}}\wedge\cdots
\wedge dx_{6}.
$$
\def\Fil{\text{Fil}}
{\bf Calculating the Picard--Fuchs Equation.}
A cohomology class in $H^{3}(V_{\lambda})$ by (7) is represented by
a differential  form
$$\eta=\sum_{i=1}^{n} \frac{P_{i}}{Q_{1}^{i}Q_{2}^{n-i}}\Omega
$$
where $\deg(P_{i})=3(n-2)$ and $n\ge 2$.
Relations among forms of this type arise from consideration
of forms $d\phi$, where
$$
\phi=\frac{\sum (x_{i}A_{j}-A_{i}x_{j}) dx_{1}\cdots
 \hat{dx_{i}}\cdots\hat{dx_{j}}\cdots dx_{6}}{Q_{1}^{i}Q_{2}^{j}}.
$$
The relations take the form
$$
\frac{
  i\sum A_{i}
	\frac{\partial Q_{1}}{\partial x_{i}}}
  {Q_{1}^{i+1}Q_{2}^{j}}
+
\frac{
  j\sum A_{i}
	\frac{\partial Q_{2}}{\partial x_{j}}}
  {Q_{1}^{i}Q_{2}^{j+1}}
\equiv
\frac{\sum\frac{\partial A_{i}}{\partial x_{i}}}{Q_{1}^{i}Q_{2}^{j}}
\pmod{\text{exact}}\eqno (9)
$$
In addition, a form with poles along only one of the forms
$Q_{i}$ is equivalent to zero:
$$
\frac{P}{Q_{i}^{j}}\Omega\equiv 0\pmod{\text{exact}}\eqno (10)
$$

We will now describe a procedure for finding canonical representations
for meromorphic forms modulo the relations (1) and (2), by constructing
an explicit representation of these relations.

Let $J_{1}$ and $J_{2}$ represent the rows of the jacobian matrix
of $(Q_{1},Q_{2})$:
$$
J_{i}=\left(\matrix
\frac{\partial Q_{i}}{\partial x_{1}} & \ldots & \frac{\partial Q_{i}}{\partial
 x_{6}}
\endmatrix
\right)
$$
If $n>2$ is an integer, we construct an $(n-1)\times 6(n-2)$
matrix $B_{n}$ as follows:
$$
B_{n}=\left(\matrix (n-2)J_{1} & 0	 & 0	     &\ldots & 0 & 0\cr
		         J_{2} & (n-3) J_{1}& 0	     & \ldots& 0 & 0\cr
			 0     &  2 J_{2} & (n-4)J_{1} & \ldots& 0 & 0\cr
			 0     &  0 	 & 3J_{2}    & \ldots & 0 & 0\cr
	    		 \vdots & \vdots & \ddots    & \ddots & 2J_{1} & 0\cr
			0	&	0 & 0	& 0 & (n-3)J_{2} & J_{1}\cr
			0	& 0	  & 0   & 0 &  0	& (n-2)J_{2}\cr
\endmatrix\right)
$$
Let $I_{n-1}$ denote the $(n-1)\times (n-1)$ identity matrix.  We
consider
the module presented by the  $(n-1)\times (8n-14)$ matrix
$K_{n}=(B_{n}\  Q_{1}I_{n-1}\ Q_{2}I_{n-1})$:
$$
S^{8n-14}\to S^{n-1}\to M^{*}_{n}\to 0,
$$
where $S$ is the graded polynomial ring in the variables $x_{1},\ldots,x_{6}$.

When $n>2$ Let $M_{n}$ denote the part of $M^{*}_{n}$ which is homogeneous of
degree $3(n-2)$, and let $M_{2}={\bold C}$.

To see the relationship between $M_{n}$ and $H^{3}(V_{\lambda})$,
suppose that $\omega$ belongs to $\Fil^{i}H^{3}(V_{\lambda})$ i.e.
 the Hodge filtration of $\omega$ is $i$.
We may represent $\omega$ in the form
$$
\omega=(\sum_{k=1}^{i+1} \frac{p_{k}}{Q_{1}^{k}Q_{2}^{i+1-k}})\Omega.
$$
Define a homogeneous map
$\phi_{i}$ from $\Fil^{i}H^{3}(V_{\lambda})$ to $S^{i+1}$
by setting $\phi_{i}(\omega)=(p_{1},\ldots,p_{i+1})$, and we
let $\overline{\phi}_{i}$ denote the composition of $\phi_{i}$ with the
projection map to $M^{*}_{i+2}$. It  follows  from
our description of the relations (9) and (10) that
  $$0 \rightarrow Fil^{i-1}H^3(V_{\lambda}) \rightarrow
 Fil^iH^3(V_{\lambda}) \rightarrow M^*_{i+2} \rightarrow 0 \eqno (11)$$
 is exact.
We may now briefly describe an algorithm for putting a cohomology
class $\omega$, presented as above, into a standard form. This standard
form will consist of elements $m_{k}(\omega)\in M_{k}$ for
$k=2,\ldots,n+2$ with the property that two forms $\omega$ and $\omega'$
represent the same cohomology class if and only if
 $m_{k}(\omega)=m_{k}(\omega')$
for all $k$ in this range.
\smallskip
{\bf Step 1.}  Compute Grobner bases for the
modules $M_{n}$ for $n=0,\ldots,i+1$.  Such a calculation provides
a canonical form for elements of $M_{i+2}$ represented as vectors in
$S^{i+1}$.
\smallskip
{\bf Step 2.} Reduce the vector $(p_{1},\ldots,p_{i+1})$ to canonical
form modulo the image of $I_{i+2}$
using Step~1.  Suppose that $m_{i+2}(\omega)$ is this canonical form.
In the reduction process, compute a vector $A$ so that
$$
\pmatrix p_{1} \cr \vdots \cr p_{i+1}\endpmatrix = m_{i+2}(\omega)+K_{i+2}A.
$$
{\bf Step 3.}  Let $A_{i,j}$ denote the subvector of $A$  consisting
of the entries $A_{i},\ldots,A_{j}$.  We denote by $\nabla\cdot A_{i,i+5}$
the usual ``divergence'' of the $6$--vector $A_{i,i+5}$ relative
to the $x_{i}$: $\nabla \cdot A_{i,i+5}=\Sigma_k {{\partial A_k}
 \over {\partial x_k}}$.
 Construct a new vector $p'$ in $S^{i}$
representing $\omega-m_{i+2}(\omega)$
(which, by the lemma, belongs to $\Fil^{i-1}$) by defining:
$$
\align
p_{1} &= \nabla\cdot A_{1,6}+A_{6i+1}\cr
p_{2} &= \nabla\cdot A_{7,12} + A_{6i+2} + A_{7i+3}\cr
      &\enspace\vdots\cr
p_{i-1} &= \nabla\cdot A_{6i-11,6i-6} + A_{7i-1} + A_{8i+1}\cr
p_{i} &= \nabla\cdot A_{6i-5,6i} + A_{8i+2}\cr
\endalign
$$
Repeat Steps~2 and~3 for $p'$, and continue
decreasing $i$ by one each time, until $i=2$.

We must apply this algorithm in one concrete situation, which we now
describe.
Define 3-forms $\omega_{i}$, for $i=2,3,\ldots$ by the formula
$$
\omega_{n}=
(-1)^{n}(n-2)!\sum_{i=1}^{n-1}
\frac{\lambda^{n}(x_{1}x_{2}x_{3})^{i-1}(x_{4}x_{5}x_{6})^{n-i-1}}{Q_{1}^{i}Q_{2
 }^{n-i}}\Omega.
$$
These define forms on the complement of ${\tilde Q}_{1} \cup {\tilde Q}_{2}$
 which, by
the residue construction, define cohomology classes on $V_{\lambda}$
invariant under the automorphism group $G_{81}$. In fact, these
forms span the space of $G_{81}$--invariant three forms on
$V_{\lambda}$, and therefore span $H^{3}(W_{\lambda})$ for the
mirror manifold.

Let $z=\lambda^{-6}$,
so that $z$ is a uniformizing parameter at $\infty$ for the
parameter space of $V_{\lambda}$.  In terms of the derivation
$$\Theta=z\frac{d}{dz}=-\frac{1}{6}\lambda\frac{d}{d\lambda} \eqno (12)$$
 we have the
following fundamental relation:
$$
\Theta\omega_{i}=-\frac{i}{6}\omega_{i}+\omega_{i+1}.\eqno (13)
$$

It follows from this relation,
$rk H^3(W(V_{\lambda})=4$, and the $G_{81}$ invariance
of the forms $\omega_{i}$ that $\omega_{6}$ is dependent on the forms
$\omega_{2},\ldots,\omega_{5}$.  By analogy with Morrison ([M2]),
we postulate a relationship of the following form:
$$
\omega_{6}(z)=\sum_{i=2}^{5}\frac{a_{i}z+b_{i}}{z-1}\omega_{i}(z)\eqno (14)
 $$  where the $a_{i}$ and $b_{i}$ are small rational numbers.
Once the $a_{i}$ and $b_{i}$ are known, it is straightforward
to compute the Picard--Fuchs equation as in [M2].
The most powerful tool available for carrying out the calculations
described in the reduction algorithm and computing
the relation (14) is the Macaulay program of
Bayer and Stillman ([Mac]).  It has one sizeable limitation
which limits its direct application to our problem --
it computes Grobner bases over a finite field, whereas at first
glance our problem requires computing over the rational function
field ${\bold C}(\lambda)$.  However, if we assume the form
of the relation we seek is as in (14), we may avoid this problem by
exploiting the Chinese Remainder Theorem:
\smallskip
{\bf Step 1.}  Set the parameter value $\lambda$ to various {\sl constant}
values $\lambda_{0}$ in the finite field ${\bold F}_{p}$.  Now use
Macaulay to apply the reduction algorithm in the corresponding fiber of the
 family
and find the relations:
$$
\omega_{6}(\lambda_{0}^{-6})=\sum
 h_{i}(\lambda_{0}^{-6})\omega_{i}(\lambda_{0}^{-6}).
$$
Here the $h_{i}$ are constants in ${\bold F}_{p}$, and
these relations are the specializations of the relation (14).
\smallskip
{\bf Step 2.} Knowledge of the values of the $h_{i}$ for, say, three
distinct $\lambda_{0}$ determines the $a_{i}$ and $b_{i}$ mod $p$.
Now repeat the calculation in Step~1 for various different choices
of $p$ (again using Macaulay), then apply the Chinese remainder theorem.
(This is not totally straightforward, since the $a_{i}$ and $b_{i}$
are rational numbers, not integers, and we have no proved {\sl a priori}
estimate on their denominators; we guessed that the denominators involved
powers of two and three, found some reasonable $a_{i}$ and $b_{i}$,
then verified that those coefficients worked for many choices of prime
$p$.)
\smallskip
Using this method, we found the following relation:
$$
\omega_{6}=\frac{z-7}{3(z-1)}\omega_{5}+
\frac{z+55}{36(z-1)}\omega_{4}+\frac{z-65}{216(z-1)}\omega_{3}+
\frac{1}{81(z-1)}\omega_{2}.  \eqno (15)$$
The associated Picard--Fuchs equation, calculated using this relation and
(13), is
the generalized hypergeometric
equation:
$$
(\Theta^{4}-z(\Theta+1/3)^{2}(\Theta+2/3)^{2})F=0\eqno (16)
$$
In particular, this implies that the monodromy at $\lambda=\infty$ is maximally
unipotent.

{\bf Computing the Yukawa Coupling.}  To determine the expansion
of the Yukawa coupling from the equation, we again follow [M2].
 The holomorphic solution $F_{0}$ to (5) is
$$
F_{0}(z)=
\sum_{n=0}^{\infty}
\left(\frac{(3n)!}{(n!)^3}\right)^2 \left(\frac{z}{3^6}\right)^{n}.
$$

We let $F_{1}$ denote the unique solution to (16) which involves
$\log(z)$ (but no higher powers of $\log(z)$) and such that
$$
s(z)=F_{1}(z)/F_{0}(z)
$$
has the property
$$
s(z)\sim \log(3^{-6}z)=-6\log(3\lambda)\quad\text{as $z\to 0$}.\eqno (17)
$$
(This is the asymptotic normalization mentioned at the beginning of the paper,
in this special case.)
If we let
$$
W=F_{0}\Theta F_{1}-F_{1}\Theta F_{0}
$$
then the Yukawa potential
 $\kappa_{sss}$,
  expressed in the canonical parameter  $q(z)=\exp(s(z))$,
 and normalized
so that its leading term is $9$ (=the degree of our Calabi--Yau family
$V_{\lambda})$
is
$$
\kappa_{sss}=-9\frac{F_{0}^4}{W^3(s(q)-1)}.
$$

To determine the predicted number of rational curves of given degree,
we write $\kappa_{sss}$ in the form
$$
\kappa_{sss}=9+\sum \frac{n_{d}d^3q^{d}}{1-q^{d}}.
$$
With our choices of normalization,
we obtain integral values for the $n_{d}$, and record them
in Table 1.

{\bf Extrapolations. }  We  know that the Picard--Fuchs
equation associated to the quintic hypersurface is the
generalized hypergeometric equation with parameters $\{1/5,\ldots,4/5\}$,
while that for the complete intersection of two cubics is the
hypergeometric equation with parameters $\{1/3,1/3,2/3,2/3\}$.  It
seems reasonable us to extrapolate from this that
the equations for the remaining
types of Calabi--Yau complete intersections are hypergeometric as well;
with parameters as given in the following table:
\smallskip
\centerline{
\vbox{\offinterlineskip\hrule
\halign{&\vrule#& \strut\quad\hfill#\hfill\quad\cr
height4pt & \omit&&\omit&&\omit&\cr
& Description && Parameters &\cr
\noalign{\hrule}
height6pt &\omit &&\omit&\cr
& $4$ quadrics in ${\bold P}^{7}$ &&
 $\{\frac{1}{2},\frac{1}{2},\frac{1}{2},\frac{1}{2}\}$&\cr
height6pt &\omit &&\omit&\cr
& $2$ quadrics and cubic in ${\bold P}^{6}$ &&
 $\{\frac{1}{2},\frac{1}{2},\frac{1}{3},\frac{2}{3}\}$&\cr
height6pt &\omit &&\omit&\cr
& $2$ cubics in ${\bold P}^{5}$ &&
 $\{\frac{1}{3},\frac{2}{3},\frac{1}{3},\frac{2}{3}\}$&\cr
height6pt &\omit &&\omit&\cr
& Quartic and quadric in ${\bold P}^{5}$ &&
 $\{\frac{1}{4},\frac{1}{2},\frac{3}{4},\frac{1}{2}\}$&\cr
height6pt &\omit &&\omit&\cr
& Quintic in ${\bold P}^{4}$ &&
 $\{\frac{1}{5},\frac{2}{5},\frac{3}{5},\frac{4}{5}\}$&\cr
height6pt &\omit &&\omit&\cr}\hrule}}

Based on this hypothesis we calculated the Yukawa potential in each of
these cases.  There are two constants which must
be chosen for each such calculation;
 one of these forces $\kappa_{sss}$ to have initial
term the degree of the variety, while the other determines the
asymptotic behavior of the coordinate $s$ in terms of the
``hypergeometric'' variable $z$ as in equation (16). In each case,
we made the choice
$$
s(z)\sim \log(z)-\sum d_{i}\log(d_{i})\eqno (18)
$$
where the $d_{i}$ are the degrees of the hypersurfaces defining the complete
intersection.
With these choices, we
obtained the correct values for the number
of straight lines in each case, and integral values for the
predicted number of rational curves.  The results of our calculations
are summarized in the Table 1.
\vfill\eject
\centerline{{\bf Table 1.}}
\centerline{{\bf Numerical Results}}
\bigskip
\centerline{Predicted Number of Rational Curves of Given Degree}
\centerline{For Various Types of Complete Intersection Calabi--Yau Manifolds}
\bigskip
\centerline{
\vbox{\offinterlineskip\hrule
\halign{&\vrule#& \strut\quad\hfill#\hfill\quad\cr
height4pt
&\omit   		&& \omit	  		&&\omit&\cr
& Degree 		&& $V_{3,3}\subset{\bold P}^{5}$ && $V_{2,4}\subset{\bold
P}^{5}$&\cr
height4pt
&\omit        		&&		\omit		 &&\omit				&\cr
\noalign{\hrule}
height6pt
&\omit 	 		&&\omit&&\omit&	  			\cr
& 1	 		&&  1053{*}	&&1280*  		&\cr
height6pt
&\omit 	 		&&\omit	  		&&\omit	&\cr
& 2	 		&& 52812	  	&&92288	&\cr
height6pt
&\omit 	 		&&\omit	  		&&\omit	&\cr
& 3	 		&& 6424326	  	&&                                    15655168	&\cr
height6pt
&\omit 	 		&&\omit	  		&&\omit	&\cr
& 4	 		&& 1139448384		&&                                   3883902528	&\cr
height6pt
&\omit 	 		&&\omit	  		&&\omit	&\cr
& 5	 		&&249787892583  	&&                                 1190923282176&\cr
height6pt
&\omit 	 		&&\omit	  		&&\omit	&\cr
& 6	 		&& 62660964509532  	&&                                417874605342336
	&\cr
height6pt
&\omit 	 		&&\omit	  		&&\omit	&\cr
& 7	 		&&17256453900822009 	&&                               160964588281789696
	&\cr
height6pt
&\omit 	 		&&\omit	  		&&\omit	&\cr
& 8	 		&& 5088842568426162960	&&
 66392895625625639488	&\cr
height6pt
&\omit 	 		&&\omit	  		&&\omit	&\cr
& 9	 		&& 1581250717976557887945&&
 28855060316616488359936	&\cr
height6pt
&\omit 	 		&&\omit	  		&&\omit	&\cr
& 10	 		&&512045241907209106828608&&
 13069047760169269024822656	&\cr
&\omit   		&&\omit	  	   	&&\omit	&\cr
\noalign{\hrule}
height 2pt
&\omit &&\omit &&\omit&\cr
\noalign{\hrule}
height4pt
& \omit&&\omit&&\omit&\cr
& Degree && $V_{2,2,2,2}\subset{\bold P}^{7}$&&$V_{2,2,3}\subset{\bold P}^{6}$&
	\cr
height 4pt
&\omit &&\omit &&\omit &\cr
\noalign{\hrule}
height6pt
 &\omit &&\omit&&\omit&   		\cr
& 1&& 512{*} &&	720*	&	\cr
height6pt
 &\omit &&\omit&&\omit&	\cr
&  2&& 9728&& 22428&		\cr
height6pt
 &\omit &&\omit&&\omit&		\cr
&  3&& 416256&&	1611504	&	\cr
height6pt
&\omit &&\omit&	&\omit&		\cr
& 4 &&25703936 &&168199200             &	\cr
height6pt
 &\omit &&\omit&&\omit&		\cr
& 5 &&  1957983744& &21676931712           &		\cr
height6pt
 &\omit &&\omit&	&\omit&	\cr
& 6 &&170535923200 &	&	3195557904564         &\cr
height6pt
 &\omit &&\omit&		&\omit&\cr
& 7 &&16300354777600 &&517064870788848       &		\cr
height6pt
 &\omit &&\omit&&\omit&		\cr
& 8 &&1668063096387072 &&89580965599606752     &	\cr
height6pt
 &\omit &&\omit&&\omit&		\cr
& 9 && 179845756064329728 &&16352303769375910848  	&\cr
height6pt
 &\omit &&\omit&	&\omit&	\cr
& 10 &&20206497983891554816 &&3110686153486233022944&	\cr
height6pt
 &\omit &&\omit&&\omit&\cr}\hrule}}
\bigskip
(*)\quad These numbers coincide with those given
 in [L] p. 52.
The number of lines on $V_{2,2,2,2}$ (resp. $V_{2,2,3}$)
 is not given explicitly there (only as part of theorem 3). It is
 easy to check that the lines  belonging to a quadric in ${\bold P}^7$ form
 a  cycle on the Grassmanian $Gr(1,7)$ of lines in ${\bold P}^7$ which is
  homologous to $4 \Omega_{4,6}$ ($\Omega_{p,q}$ denotes the Schubert
 cycle consisting of lines in a generic ${\bold P}^q$ intersecting
 generic ${\bold P}^p \subset {\bold P}^q$ ).
 Its 4-fold self--intersection equals $512$, which gives the number
 of lines on $V_{2,2,2,2}$. On the other hand, the lines in ${\bold P}^6$
 which belong to a generic hypersurface of degree $3$ (resp. $2$) form
 the cycle in $Gr(1,6)$ homologous to $18 \Omega_{2,5}+27 \Omega_{3,4}$
(resp. $4 \Omega_{3,5}$). The intersection index: $(18 \Omega_{2,5}+
 27 \Omega_{3,4})(4\Omega_{3,5})^2$ equals $720$ which gives the
 number of lines on $V_{3,2,2}$.
\bigskip
\centerline {\bf References}
\bigskip
\par [COGP] P.Candelas, X. de la Ossa, P.Green, L.Parkes,
 Nuclear Physics. B359 (1990). 21.
\par [CL] P.Candelas, M.Lynker, R.Schimmrigk,
Calabi--Yau  manifolds in weighted
 ${\bold P^4}$. Nuclear Physics, B341 (1990), 383-402.
\par [DHVW] L.Dixon, J.Harvey, C.Vafa, E.Witten, Strings on
 orbifolds, I, Nuclear Phys. B261, 678-686 (1985) and II (ibid)
 B 274, 285-314 (1986).
\par [G] P.Griffiths, Periods of integrals on algebraic manifolds:
 summary of results and discussions of open problems, Bull. AMS.
 1970, 76, 228-296.
\par [GP] B.Greene and M.Plesser, Duality in Calabi Yau moduli space
, Nuclear Physics  B338 (1990).
\par [HH] F.Hirzebruch and T.Hofer, On Euler number of an orbifold,
 Math. Ann. 286, 255-260, (1990).
\par [L] A.Libgober, Numerical characteristics of systems of
 straight lines on complete intersections, Math. Notes, 13(1973) p. 51--56,
 Plenum
Publishing, translated from Math. Zametki 13(1973) p. 87--96.
\par [Mac] D. Bayer and M. Stillman, {\it Macaulay: A system
for computation in algebraic geometry and commutative algebra}.
Source and object code available for Unix and Macintosh computers.
Contact the authors, or download from {\tt zariski.harvard.edu} via
anonymous ftp.
\par [M1] D.Morrison, Mirror symmetry and rational curves on
 quintic threefolds: a guide for mathematicians,
Preprint, Duke University, 1991.
\par [M2] D.Morrison, Picard Fuchs equations and the mirror map for
 hypersurfaces, Preprint, Duke University, 1991.
\par [S] R.Schimmrigk, The construction of mirror symmetry. Preprint,
 Sept. 1992, Heidelberg and CERN.
\par P.S. Besides calculation in [L] of the number of lines on generic
 complete intersection of arbitrary dimension on the case when it is
 finite the case of lines on three dimensional complete intersections
 with $K=0$ also treated in S.Katz paper in Math. Zeit. 191 (1986)
 .293-296. We are thanking S.Katz for pointing this out.
\end